 \documentclass[draft]{article}

\usepackage{amsmath,amsfonts,amsthm,amssymb,amscd,cancel,color,mathabx}
\usepackage{enumitem}
\usepackage{verbatim}
\usepackage[dvipdfmx]{graphicx}
\usepackage{ulem,color}

\renewcommand{\Re}{\mathop{\rm Re}\nolimits}
\renewcommand{\Im}{\mathop{\rm Im}\nolimits}

\theoremstyle{plain}
\newtheorem{theorem}{Theorem}[section]
\newtheorem{lemma}[theorem]{Lemma}
\newtheorem{proposition}[theorem]{Proposition}

\theoremstyle{definition}

\theoremstyle{remark}
\newtheorem{remark}[theorem]{Remark}
\newtheorem{notation}[theorem]{Notation}

\newcommand{\R}{{\mathbb R}}

\newcommand{\Z}{{\mathbb Z}}

\def\im{{\rm i}}

\newcommand{\C}{\mathbb{C}}

\def\({\left(}
\def\){\right)}
\def\<{\left\langle}
\def\>{\right\rangle}
\newcommand{\sech}{{\mathrm{sech}}}

\newcommand{\Span}{{\mathrm{Span}}}

\makeatletter

\numberwithin{equation}{section}

\begin{document}

\title{On small energy solutions of the Nonlinear Schr\"odinger Equation in 1D with a generic trapping potential with a single eigenvalue.}

\author{Scipio Cuccagna, Masaya Maeda }
\maketitle

\begin{abstract} We prove in  dimension $d=1$ a result similar to Soffer and Weinstein  Jour. Diff. Eq. 98 (1992) capturing for pure power nonlinearities the whole range of exponents $p>1$.  The proof is based on the virial inequality of   Kowalczyk  \textit{{et al.}} J. Eur. Math. Soc. (JEMS) 24 (2022) with smoothing estimates like in Mizumachi   J. Math. Kyoto Univ. 48 (2008).
\end{abstract}

\section{Introduction}

We consider  the   nonlinear Schr\"odinger equation  (NLS) on the line
\begin{align}\label{eq:nls1}&
  \im \partial _t  u  =Hu +  f(u)    ,\quad (t,x)\in \R\times \R ,\text{with $u(0)=u_0 \in H^1 (\R , \C )$   and $  f(u )= g(|u|^2)  u$}
\end{align}
where $H=-\partial ^2_x +V$   is a Schr\"odinger operator with $V\in \mathcal{C}^1(\R, \R )$, that is real valued, with $|V^{(k)}(x) |\lesssim e^{-\kappa _0 |x|}$ for $k=0,1$ and some $\kappa _0>0$. We will assume that 0 is not a resonance for $H$, that is
\begin{align}
  \label{eq:nonres}   Hu=0    \text {  with }   u\in L^\infty (\R   ) \Longrightarrow   u=0.
\end{align}
 It is well known that then $H$ has a finite number of eigenvalues and that they are of dimension one and that for the continuous spectrum we have
 $\sigma _e(H) =[0,+\infty )$. We will assume that $H$ has exactly one single eigenvalue, which we denote by $-\lambda$ and we will assume $-\lambda <0$.
 We will denote by $\varphi \in \ker (H +\lambda )$ a corresponding eigenfunction with $\| \varphi \| _{L^2}=1$,  $\varphi >0$  everywhere, with $\varphi (x)\lesssim e^{-|x|\sqrt{\lambda}}$,  \cite  {DT}.
  We will consider the spectral decomposition
  \begin{align}
    \label{eq:specdec} L^2(\R, \C  ) = \Span \{ \varphi  \}\oplus L_c(H)
  \end{align}
 denoting by $P u: = \<u,\varphi\>\varphi +\<u,\im \varphi\>\im\varphi$ the corresponding projection on $\Span \{ \varphi  \}$ and $P_c=1-P$
 where   \begin{equation}\label{eq:bilf}
  \( u , v \) =   \int _{\R  } u(x) \overline{v}(x) dx      \text{ and }  \langle u , v \rangle = \Re \( u , v \) \text{  for  $u,v:\R \to \C$ }.
 \end{equation}
 We will also use the following notation.
\begin{itemize}
  \item Given a Banach space $X$, $v\in X$ and $\varepsilon>0$ we set
$
D_X(v,\varepsilon):=\{ x\in X\ |\ \|v-x\|_X<\varepsilon \}.
$
\item
For $\gamma\in \R$ we set
\begin{align}
L^2_\gamma&:=\{u\in \mathcal S'(\R,\C)\ |\ \|u\|_{L^2_\gamma}:=\|e^{\gamma|x|}u\|_{L^2}<\infty\},\label{L2g}\\
H^1_\gamma&:=\{u\in \mathcal S'(\R,\C)\ |\ \|u\|_{H^1_\gamma}:=\|e^{\gamma|x|}u\|_{H^1}<\infty\}.\label{H1g}
\end{align}
\item   For $s\in \R$ we set    and   $\kappa \in (0,1)$    fixed in terms of $p$  and  small enough, we consider
\begin{align}\label{eq:l2w}&
\|  {\eta} \|_{  L ^{p,s}} :=\left \|\< x \> ^s \eta \right \|_{L^p(\R )}   \text{  where $\< x \>  := \sqrt{1+x^2}$,}\\&
\| \eta  \|_{  { \Sigma }_A} :=\left \| \sech \(\frac{2}{A} x\) \eta '\right \|_{L^2(\R )} +A^{-1}\left \|    \sech \(\frac{2}{A} x\)  \eta   \right\|_{L^2(\R )}  \text{ and} \label{eq:normA}\\& \|  {\eta} \|_{ \widetilde{\Sigma} } :=\left \| \sech \( \kappa   x\)   {\eta}\right \|_{L^2(\R )} . \label{eq:normk}
\end{align}

\end{itemize}
We assume    there exist  $p>1$  and $C>0$ such that for $k=0,1,2 $ we have
\begin{align}\label{nonlinearity}
 |g ^{(k)}(s)| \leq C |s|^{\frac{p-1}{2}-k}  \text{    for all $s\in(0,1]$.}
\end{align}
The following can be proved like in   \cite[Proposition 1.2]{CM19SIMA}.
\begin{proposition}[Bound states]\label{prop:bddst} Let $p>0$.
Then there exist $a_0>0$, $\varepsilon _0>0$ and $C>0$ such that there exists a unique $Q[\cdot ]\in C^1(D_\C(0,\varepsilon _0),H^1_{\gamma_0})$ satisfying the gauge property
\begin{align}\label{Q:gauge}
Q[e^{\im \theta}z]=e^{\im \theta}Q[z],
\end{align}
such that there exists $E\in C([0,\varepsilon _0^2),\R)$ such that
\begin{align}\label{eq:sp}
H  Q[z]+g(|Q[z]|^2)Q[z]=E(|z|^2)Q[z],
\end{align}
and for $j=1,2$,
\begin{align}\label{prop:bddst1}
\|Q[z]-z\varphi\|_{H^1_{a_0}}\leq C |z|^{p},\ \|D_jQ[z]-\im^{j-1}\varphi\|_{H^1_{a_0}}\leq C |z|^{ p-1  },\
\left|E(|z|^2)+\lambda \right|\leq C |z|^{ p-1  }.
\end{align}
\end{proposition}
\qed

\noindent  Energy $\mathbf{E}$ and Mass $\mathbf{Q}$ are invariants of \eqref{eq:nls1}, where
  \begin{align}\label{eq:energy}
& \mathbf{E}( {u})=\frac{1}{2}   \< Hu ,u  \> +\int_{\R}  G (|u|^2)    \,dx   \text {  where }  G(0) =0 \text{ and }     G '(s) \equiv g(s) , \\&   \label{eq:mass} \mathbf{Q}( {u})=\frac{1}{2}   \| u  \| ^2 _{L^2\( \R \) }.
\end{align}
In this paper we prove the following result.  Near the origin  $H^1(\R ) $     problem  \eqref{eq:nls1}  is globally well posed and there exist constants $\delta _0>0$ and $C_0>0$ such that for
$ \| u _0\| _{H^1}<\delta  _0 $  then $ \| u (t)\| _{H^1} \le C_0 \| u _0\| _{H^1}$  for all times, see Cazenave \cite{CazSemi}.

\begin{theorem} \label{thm:asstab}   There exists an $a>0$ such that   for any $\epsilon  >0$
there exist  $\delta   >0$   such that  for $ \| u _0\| _{H^1}<\delta   $ the  solution  $u(t)$ of  \eqref{eq:nls1} can be written uniquely  for all times as
 \begin{equation}\label{eq:small en1}
\begin{aligned}&    u(t)= Q[z (t)]+\eta (t) \text{ with $\eta (t) \in
P_c H^1$,}
\end{aligned}
\end{equation}
such that  we have
\begin{align}&    \label{eq:small en2} |z(t)|+\| \eta (t)\|  _{H^1}\le \epsilon  \text{ for all $t\in [0,\infty ) $, }
 \\ & \label{eq:asstab2}
       \int _0 ^{+\infty }   \|    e ^{-a \< x \>}\eta \|_{H^1 }^2 dt \le   \epsilon ^2  \\&   \label{eq:asstab20}  \lim _{t\to +\infty  }  \|  e^{- a \< x\>}   \eta (t ) \| _{L^2(\R )}     =0
\end{align}
 and     there exists  a $r  _+\ge 0$ such that \begin{equation}\label{eq:small en22} \begin{aligned}&       \lim _{t\to +\infty } |z (t) |=r _+. \end{aligned} \end{equation}
\end{theorem}

\begin{remark}\label{rem:asstab3--}
Theorem \ref{thm:asstab}  implies   that for $\delta $ sufficiently small    there exists a function $\vartheta :\R \to \R$ such that
\begin{align}\label{eq:rem:asstab3--1}
  \lim _{t\to +\infty }  u(t) e^{ -\im \vartheta (t)   }=   Q[r _+]  \text{  in }L^\infty _{\text{loc}}(\R ) .
\end{align}
\end{remark}

 A version of Theorem \ref{thm:asstab} for $p>2$  but where the nonlinearity  is  multiplied by an appropriately decaying function in $x\to \infty $ in Weder \cite{weder}.  The case of the cubic nonlinearity was considered by Gong Chen \cite{Chen} and by Masaki et al. \cite{MMS1} if $V(x)= - c_0 \delta (x)$, where $c_0>0$ and $\delta (x)$ the Dirac delta concentrated in 0. We considered the     cubic nonlinearity in \cite{CM2109.08108} with a more general set of eigenvalues than the one considered here.  In    \cite{CM2109.08108} the nonlinearity could be different but needs to be sufficiently regular and the regularity depends on the ratio of the distances of the first excited state from the ground state and from 0. Notice that in  \cite[Assumption 1.13]{CM2109.08108} there is an additional repulsivity hypothesis on a potential obtained after Darboux transformations, which is actually unnecessary,  in fact is not required if we replace the second virial inequality in \cite{CM2109.08108} with smoothing estimate like the ones considered here.

\noindent Here we will focus only on the more delicate  case $1<p\le 2$, where the nonlinearity is stronger. In \cite{CM19SIMA} we considered a special case  of Theorem \ref{thm:asstab} with $V(x)= - c_0 \delta (x)$. We improve that result because  in  \cite{CM19SIMA} we proved a version of
\eqref{eq:small en22} only for $p>2$.

Equation \eqref{eq:nls1} was considered first for dimension $d=3$ by       Soffer  and Weinstein \cite{SW1,SW2}  and  Pillet  and Wayne  \cite{pillet}.    Of the subsequent literature up until 2020 we highlight
Gustafson et  al.   \cite{GNT} for the case $d=3$ and the adaption of \cite{GNT} for $d=1,2$ by Mizumachi \cite{mizu08}, while we refer
 for a  long list of references on the asymptotic stability of ground states of the NLS   up until 2020 to   our  survey \cite{CM21DCDS}.  Various papers have been written recently proving dispersive estimates of small solutions in the case when $H$ has no eigenvalues, see \cite{Delort}--\cite{GoPu1}, \cite{LP0}--\cite{MMS2} and \cite{naumkin2016}. The  papers proving  dispersive estimates (for a smaller class of solutions), which appear harder to prove than Theorem \ref{thm:asstab},   require a certain degree of regularity for the nonlinearity. Notice that in our $H^1\( \R \)$ framework all results  must be invariant by translation in time, a condition satisfied by Theorem \ref{thm:asstab},  and so it is not possible to have a fixed  rate  of time  decay    $\| \eta  (t)  \| _{L^\infty (\R )}$  for all our solutions.
Li and  Luhrmann \cite{LiLu} partially recover  the result in  \cite{cupe2014} with methods which do not exploit explicitly     the integrable structure of the cubic NLS.

 Like in  \cite{CM24D1} our proof is based on the ideas of Kowalczyk  \textit{{et al.}}  \cite{KM22}--\cite{KMMvdB21AnnPDE}  where dispersion is proved by means of two distinct virial estimates. Like in
\cite{CM24D1}  here we perform just the first, high energy, virial estimate while we substitute the second, low energy, virial  estimate by a Kato smoothing estimate. Very important to us has been the treatement of  Kato smoothing by  Mizumachi \cite{mizu08}. The methods of   Kowalczyk  and our particular use of  smoothing estimates do not require much regularity of the nonlinearity and allow us to consider equation \eqref{eq:nls1} for any $p>1$.
   For an earlier example of this for a different model, using both virial inequalities,  see \cite{CMMS23}.  Our methods allow us to exploit the fact that, under our hypotheses on the Schr\"odinger operator $H$, at small energy radiation disperses to infinity more than small energy solutions of the constant coefficients linear Schr\"odinger equation.  The  two classes of solutions are different, see \cite{MN,strauss}.  Notice that  for $H=- \partial _x^2$ and  $p\ge 3$ by Ozawa \cite{ozawa} for any $\eta _+\in L^{2,2}(\R )$, the space defined in  \eqref{eq:l2w}, and both for focusing and defocusing equation   there exists a solution $u$ of \eqref{eq:nls1} with $\| \eta(t) - e^{\im S(t,\cdot )}e^{\im t\partial ^2_x}\eta _+\| _{L^2(\R )}\le O(t^{-\alpha}) $ with $\alpha >1/2$ and for an appropriate real phase $S(t,x )$, which is  equal to 0 if $p>3$. Since   for  $|x|\le 1$   we have $e^{\im t\partial ^2_x}\eta _+( x) = \frac{1}{ \sqrt{2\im t}  } e^{\im \frac{x^2}{4t}} \widehat{\eta}_+(0) +O(t^{-\alpha}) $ with $\alpha >1/2$, in the generic case with   $\widehat{\eta}_+(0)\neq 0$  we conclude that $\|   \eta  \| _{L^2 (\R _+ , L^2(|x|\le 1))} =+\infty$. The fact that here  we have  \eqref{eq:asstab2} is ultimately due to the fact that
our $H$ is generic and so is better behaved in terms of the Kato smoothing than $-\partial _x^2$   which is non generic has a resonance  at 0.  The second virial inequality in various of the papers of  Kowalczyk  \textit{{et al.}}  and in \cite{CMMS23} is also due to the fact a certain linearized operators are generic  in the framework  of these papers.

Turning on open problems, we do not know how to prove Theorem \ref{thm:asstab}   for low $p$'s in dimension $d\ge 2$.  Another interesting problem would be   improve      \cite{Cu} and  prove a version of Theorem \ref{thm:asstab} on the line for a potential $V=q_1+q_2$ for $q_1$ periodic and $q_2$ a small and rapidly decreasing to 0 at infinity. Unfortunately  our proof relies on the fact that the background is flat and $V$ decays rapidly, whether exponentially like here or at a sufficiently rapid polynomial decay.

\begin{notation}\label{not:notation} We will use the following miscellanea of  notations and definitions.
\begin{enumerate}


\item Like in the theory of     Kowalczyk et al.   \cite{KMM3},      we   consider constants  $A, B,  \epsilon , \delta  >0$ satisfying
 \begin{align}\label{eq:relABg}
\log(\delta ^{-1})\gg\log(\epsilon ^{-1}) \gg     A  \gg    B^2\gg B  \gg 1.
 \end{align}
 Here we will take $ A\sim B^3$, see Sect. \ref{sec:smooth1} below,  but in fact  $ A\sim B^n$ for  any   $n>2$ would make no difference.

 \item The notation    $o_{\varepsilon}(1)$  means a constant   with a parameter $\varepsilon$ such that
 \begin{align}\label{eq:smallo}
 \text{ $o_{\varepsilon}(1) \xrightarrow {\varepsilon  \to 0^+   }0.$}
 \end{align}

\item Given two Banach spaces $X$ and $Y$ we denote by $\mathcal{L}(X,Y)$ the space of continuous linear operators from $X$ to $Y$. We write $\mathcal{L}(X ):=\mathcal{L}(X,X)$.

\item We have the following elementary formulas,
\begin{align} \label{eq:derf1}&
  Df(u)X=\left . \frac{d}{dt}  f\( u+tX \)  \right | _{t=0}= g\( |u  | ^{2}\) X + 2  g '\( |u  | ^{2}\) u \< u, X\> _\C  .
\end{align}

\item Following  the framework in Kowalczyk et al. \cite{KMM3} we   fix an even function $\chi\in C_c^\infty(\R , [0,1])$ satisfying
\begin{align}  \label{eq:chi} \text{$1_{[-1,1]}\leq \chi \leq 1_{[-2,2]}$ and $x\chi'(x)\leq 0$ and set $\chi_C:=\chi(\cdot/C)$  for a  $C>0$}.
\end{align}

\item We set
\begin{align}\label{eq:C+-}
 \C _\pm := \{ z\in \C : \pm \Im z >0   \} .
\end{align}

 \end{enumerate}

\end{notation} \qed

\section{Notation and coordinates} \label{section:set up}

 We have the following  ansatz, which is an elementary consequence of the Implicit Function Theorem, see \cite[Lemma 2.1]{CM19SIMA}.
\begin{lemma}\label{lem:decomposition}
There exist $c _0 >0$ and $C>0$ such that for all $u \in H^1$   with $\|u\|_{H^1}<c  _0 $, there exists a unique  pair $(z,\eta )\in   \C \times  P_c H^1 $
such that
\begin{equation}\label{eq:ansatz}
\begin{aligned} &
u= Q[z]+\eta  \text{ with }  |z |+\|\eta  \|_{H^1}\le C \|u\|_{H^1} .
\end{aligned}
\end{equation}
 The map  $u \to (z,\eta )$  is in $C^1 (D _{H^1}(0, c _0 ),
 \C   \times  H^1 )$.
\end{lemma}
 \qed

The proof of Theorem \ref{thm:asstab} is mainly  based on the following continuation argument.

\begin{proposition}\label{prop:continuation}
There exists  a    $\delta _0= \delta _0(\epsilon )   $ s.t.\  if
\begin{align}
  \label{eq:main2} \| \eta \| _{L^2(I, \Sigma _A )} +  \| \eta \| _{L^2(I, \widetilde{\Sigma}  )} + \| \dot z +\im E z \| _{L^2(I  )} \le \epsilon
\end{align}
holds  for $I=[0,T]$ for some $T>0$ and for $\delta  \in (0, \delta _0)$
then in fact for $I=[0,T]$    inequality   \eqref{eq:main2} holds   for   $\epsilon$ replaced by $   o _{\epsilon }(1) \epsilon $.
\end{proposition}
Notice that this implies that in fact the result is true for $I=\R _+$.   We will split the proof of Proposition \ref{prop:continuation} in a number of partial results obtained assuming the hypotheses of Proposition  \ref{prop:continuation}.
\begin{proposition}\label{prop:modpar} We have
  \begin{align}&
  \label{eq:modpar1} \| \dot  z   +\im E (|z|^2) z \|  _{L^2(I  )}  \lesssim  \delta^{p-1} \epsilon   ,\\&    \label{eq:modpar3}   \|\dot z \|  _{L^\infty(I  )} \lesssim   {\delta} .
\end{align}

\end{proposition}

\begin{proposition}[Virial Inequality]\label{prop:1virial}
 We have
 \begin{align}\label{eq:sec:1virial1}
   \| \eta \| _{L^2(I, \Sigma _A )} \lesssim  A \delta  + \| \dot  z   +\im E (|z|^2) z \|  _{L^2(I  )}   + \| \eta \| _{L^2(I, \widetilde{\Sigma}   )} + \epsilon ^2  .
 \end{align}
 \end{proposition}
Notice that $A \delta \ll  A ^{-1}\epsilon ^2 =o _{B^{-1}} (1)  \epsilon ^2 $  in \eqref{eq:sec:1virial1}.

\begin{proposition}[Smoothing Inequality]\label{prop:smooth11}
 We have
 \begin{align}\label{eq:sec:smooth11}
   \| \eta \| _{L^2(I, \widetilde{\Sigma}   )} \lesssim  o _{B^{-1}} (1)  \epsilon    .
 \end{align}
 \end{proposition}

%

\textit{Proof of Theorem \ref{thm:asstab}.}
It is straightforward that Propositions \ref{prop:modpar}--\ref{prop:smooth11}  imply Proposition \ref{prop:continuation} and thus the fact that we can take $I=\R _+$ in all the above inequalities. This in particular implies \eqref{eq:asstab2}.

\noindent We next focus on the limit \eqref{eq:asstab20}. We first rewrite our equation,
entering the  ansatz \eqref{eq:ansatz} in \eqref{eq:nls1} and using   \eqref{eq:sp}    we obtain
\begin{align}       &
    \im \dot \eta   +\im    D_{z} Q [z ]  \(\dot  z   +\im E z \)     = H   \eta  +   \( f(  Q [z ]  + \eta  )  -  f(  Q [z ]  )  \)     .   \label{eq:nls3}
\end{align}
Then we can proceed like in \cite{CM24D1} considering
  \begin{equation*}
  \begin{aligned}
   \mathbf{a}(t) &:= 2 ^{-1}\|  e^{- \gamma \< x\>}   \eta (t)  \| _{L^2(\R )} ^2
  \end{aligned}
  \end{equation*}
and  by  obvious cancellations and       orbital stability getting  \small
\begin{align}\label{eq:intparts}
  \dot {\mathbf{a}}  &=    \< \left [   e^{- 2\gamma \< x\>} , \im  \partial _x ^2  \right ] \eta , \eta \>
   -\<    D_{z} Q [z ]   \(\dot  z+  \im E z \)  , e^{- 2a\< x\>} \eta    \>   \\& + \<
  e^{- \gamma \< x\>}        \im  \( f( Q [z ]   + \eta  )  -  f( Q [z ]    )  \)  , e^{- \gamma\< x\>} \eta    \> =O( \delta  ^2 ) \text{ for all times}.\nonumber
\end{align}\normalsize
   Since we already know from \eqref{eq:asstab2} that $\mathbf{a}\in L^1(\R )$, we conclude that
$ \mathbf{a}(t) \xrightarrow{t\to +\infty} 0$.
Notice that the integration by parts in \eqref{eq:intparts} can be made rigorous   considering that
if $u_0\in H^2(\R )$ by the well known regularity result by Kato, see \cite{CazSemi}, we have $\eta \in C^0\( \R , H^2 (\R )\)$ and the above argument is correct and by a standard  density argument the result  can be extended to $u_0\in H^1(\R )$.

We prove \eqref{eq:small en22}.
Here, notice that by orbital stability we can  take   $a >0$ such that we have the  following, which will be used below,
\begin{equation}\label{eq:uniformity1}
  e^{-2a\<x\>} |z| \gtrsim    \max \{ |Q [z ]| , |Q [z ]| ^{2p-1} \}  \text {  for all } t\in \R .
\end{equation}
Since $\mathbf{Q}(Q [z  ]) $ is strictly  monotonic in $|z|$, it suffices to show the  $\mathbf{Q}(Q [z (t)])$ converges as $t\to\infty$.
From the conservation of $\mathbf{Q}$, the exponential decay of $\phi[\omega,v,z]$, \eqref{eq:asstab20} and \eqref{eq:asstab2}, we have
\begin{align}\label{eq:equipart1}
\lim_{t\to \infty}\(\mathbf{Q}(u_0)-\mathbf{Q}(Q [z (t)])-\mathbf{Q}(\eta (t))\)=0.
\end{align}
Thus, our task is now to prove $\frac{d}{dt}\mathbf{Q}(\eta) \in L^1$, which is sufficient to show the convergence of $\mathbf{Q}(\eta)$.
Now, from \eqref{eq:nls3}, we have
\begin{align*}
\frac{d}{dt}\mathbf{Q}(\eta)&=\<\eta,\dot{\eta}\>
=-\<\eta, D_{z} Q [z ]   \(\dot  z+  \im E z \)      \>    +\<\eta,  \im  \( f( Q [z ]   + \eta  )  -  f( Q [z ]    )  \)\> \\&
=I+II  .
\end{align*}
By the bound of the 1st and the 3rd term of \eqref{eq:main2}, we have $I\in L^1(\R _+)$.
  To show  $II\in L^1(\R _+)$ we partition  for $s\in (0,1)$ the line   where $x$ lives as
\begin{align*}&
   \Omega_{1,t ,s}=\{x\in \R\ |\  | s \eta(t,x)|\leq 2| Q[ z(t)]|\}  \text{ and }\\&  \Omega_{2,t,s }=\R\setminus \Omega_{1,t,s }=\{x\in \R\ |\ | s\eta(t,x)|> 2|Q[ z(t)] |\} ,
\end{align*}
Then, we have
\begin{align*}
II(t)&=\sum_{j=1,2}   \int_0^1ds \int_{\Omega_{j,t ,s }} \< \im  D  f( Q[ z(t)]+ s\eta (t))  \eta (t) , \eta (t)        \> _\C     \,dx\\& =:II_1(t)+II_2(t).
\end{align*}
For $II_1$,  by  \eqref{nonlinearity} and  \eqref{eq:derf1}, we have \small
\begin{align*}
|II_1(t)|&\lesssim   \int_0^1ds\int_{\Omega_{1,t,s }}  \left |   \< \im  D  f(Q[ z(t)] + s\eta (t))  \eta (t) , \eta (t)        \> _\C        \right | dx   \\& \lesssim \int_0^1ds\int_{\Omega_{1,t ,s}}  |Q[ z(t)]   +s \eta (t) | ^{ p-1}|\eta (t)|^2\,dx \lesssim \int_{\R} |Q[ z(t)] |^{ p-1}|\eta (t)|^2\,dx \\&  \lesssim \|  e^{- \gamma \<x\>}    \eta (t)\| _{L^2 (\R )}  ^2 \in L^1(\R _+ ).
\end{align*}\normalsize
Turning to $II_2$,  by exploiting $ \<  D  f( s \eta (t)) \eta (t)           ,\im  \eta (t) \> _\C  \equiv 0  $,  which can be easily checked from   \eqref{eq:derf1},
we write \small
\begin{align*}
   II_2(t) &= - \int _{[0,1]  }     ds \int_{\Omega_{2,t,s }} \<  D  f(Q[ z(t)]+ s\eta (t)) \eta (t) -D  f( s\eta (t)) \eta (t)           ,\im  \eta (t) \> _\C      dx \\&  = - \int _{[0,1] ^2 }  d\tau ds  \int_{\Omega_{2,t,s }} \<  D  ^2 f(\tau Q[ z(t)]+s \eta (t))  (  Q[ z(t)] ,\eta (t))             ,\im  \eta (t) \> _\C      dx .
\end{align*}\normalsize
Then by  $|s \eta(t,x)|> 2|Q[ z(t)]|$,   $p\in (1,2]$, we get the following, which completes the proof of $II\in L^1(\R _+)$ and \eqref{eq:small en22},
\begin{align*}
   |II_2(t)| &\lesssim    \int _{[0,1] ^2 }  d\tau  ds \int_{\Omega_{2,t,s }}       | s \eta (t)  | ^{ p-2}  |  Q[ z(t)]  |   \ |\eta (t) |^2     dx \\& \lesssim   \int _{[0,1]   }    ds \ s^{ p-2}  \int_{\Omega_{2,t,s  }}       |  \eta (t)  | ^{ p-2}  | Q[ z(t)] | ^{ p-1}  |  \eta (t)  | ^{2- p}  \ |\eta (t) |^2     dx \\&\lesssim    \|  e^{- a\<x\>}    \eta (t)\| _{L^2 (\R )}  ^2 \in L^1(\R _+ ).
\end{align*}

\qed

\section{Proof of Proposition  \ref{prop:modpar} }\label{sec:modpar}

Proposition  \ref{prop:modpar} is an immediate consequence of \eqref{eq:main2} and the followig lemma.

\begin{lemma}\label{lem:lemdscrt} We have the following estimate,
 \begin{align}   |\dot  z   +\im E (|z|^2) z  | \lesssim    \delta ^{p-1} \|   \eta \| _{\widetilde{\Sigma}}    .   \label{eq:discrest2}
\end{align}

\end{lemma}
\proof
Applying $\< \cdot ,    \Theta \varphi \>  $  with $ {\Theta } =1, \im $   to \eqref{eq:nls3}
 we have the following, where the cancelled terms are null   by Lemma \ref{lem:decomposition},
\begin{align*}    &  \cancel{\<  \im \dot \eta  ,   \Theta \varphi   \> } +    \<   \im    D_{z} Q [z ]  \(\dot  z   +\im E z \)   , \Theta \varphi  \>
       =    \cancel{\< H\eta ,    \Theta \varphi     \>  }  \\&   -\<     f( Q [  z ] + \eta  )  -  f( Q [  z ] ) -D f( Q [  z ]   ) \eta ,  \Theta \varphi \> -  \<     D f( Q [  z ]   ) \eta ,  \Theta \varphi \>  .
   \end{align*}
 We have
 \begin{align} \label{eq:discrestnonlin1}
   | \<     D f( Q [  z ]   ) \eta ,  \Theta \varphi \>| \lesssim |z| ^{p-1}    \|   \eta \| _{\widetilde{\Sigma}}.
 \end{align}
 Next, we claim    \begin{align} \label{eq:discrestnonlin1zj}
 | \<     f( Q [  z ] + \eta  )  -  f( Q [  z ]  ) -D f( Q [  z ]  ) \eta ,      \Theta \varphi    \> |\lesssim \|  \eta  \| _{\widetilde{\Sigma}} ^ p.
\end{align}
 Set  for $s\in (0,1)$
\begin{align*}&
   \Omega_{1,t ,s}'=\{x\in \R\ |\  2| s \eta(t,x)| \le |Q [ z(t)]|\}  \text{ and }\\&  \Omega_{2,t,s }'=\R\setminus \Omega_{1,t,s } '=\{x\in \R\ |\ 2| s\eta(t,x)|>  |Q [ z(t)] |\}
\end{align*}
and  split
\begin{align*} &
  \<     f( Q [  z ]  + \eta  )  -  f(Q [  z ]   ) -D f( Q [  z ]  ) \eta ,    \Theta  \varphi \> =I_1(t)+I_2(t) \text{ for}
   \\& I_j(t)=  \int _0^1 ds \int _{\Omega_{j,t,s }'}  \<   \left [ D f( Q [  z ]  + s\eta  )    -D f( Q [  z ]  ) \right ] \eta ,   \Theta \varphi \>  _\C dx .
\end{align*}
We have
\begin{align*} &
   I_1 (t)=  \int _{[0,1]^2} d \tau \  ds  \ s\int _{\Omega_{1,t,s }'}  \<   D ^2f( Q [  z ]  + s \tau \eta  )  ( \eta , \eta ) ,    \Theta \varphi   \>  _\C dx
\end{align*}
with
\begin{align*}
  |I_1(t)|& \lesssim    \int _{[0,1]^2} d \tau \  ds  \int _{\Omega_{1,t,s }'}   \varphi  | Q [  z ] | ^{ p-2} |\eta |  ^{2-p +p} d x \lesssim \|  \eta  \| _{\widetilde{\Sigma}} ^p.
\end{align*}
We have the following, which completes the proof of  \eqref{eq:discrestnonlin1zj} and with \eqref{eq:discrestnonlin1} yields \eqref{eq:discrest2},
\begin{align*}
  |I_2(t)|& \lesssim    \int _{[0,1] }  \  ds  \int _{\Omega_{2,t,s }'} \varphi    |\eta | ^p  d x \lesssim
   \|  \eta  \| _{\widetilde{\Sigma}} ^p  .
\end{align*}

 \qed

\section{High energies: proof of Proposition \ref{prop:1virial}}\label{sec:vir}

Following  the framework in Kowalczyk et al. \cite{KMM3}  and
using the function $\chi$  in \eqref{eq:chi}
we  consider the   function
\begin{align}\label{def:zetaC}
\zeta_A(x):=\exp\(-\frac{|x|}{A}(1-\chi(x))\)   \text {  and   }  \varphi_A (x):=\int_0^x \zeta_A^2(y)\,dy
\end{align}
  and  the vector field
\begin{align}\label{def:vfSA} \  S_A:=   \varphi_A'+2 \varphi   _{ A}\partial_x .
\end{align}
Next  we introduce   \begin{align}\label{def:funct:vir}
\mathcal{I}_A := 2^{-1}\<  \im \eta    ,S_A  \eta  \>     .
\end{align}
Notice that $|\mathcal{I}_A (t) |\lesssim  A \delta ^2$ for any $t\in \R$.

\begin{lemma}\label{lem:1stV1}
There exists a fixed constant $C>0$ s.t. for an arbitrary small number
\begin{align} &
\|   \eta  \|_{\Sigma _A}^2   \le C\left [  \dot{\mathcal{I}}_A +  \|    \eta \| _{\widetilde{\Sigma}} ^2  +|\dot  z   +\im E z | ^2  \right ].   \label{eq:lem:1stV1}
\end{align}
\end{lemma}
\proof
From \eqref{eq:nls3}   we obtain
\begin{align*}
    \dot {\mathcal{I}}_A   =&-\<  \dot \eta,    \im S_A\eta \>   = -\<  \partial _x^2    \eta ,S_A\eta  \>+\<V\eta,S_A\eta\>    \\& -  \<    f( Q [z ] + \eta  )  -  f( Q [z ]    )     ,S_A\eta  \>   + O\(  |\dot  z   +\im E z|   \|    \eta \| _{\widetilde{\Sigma}}  \) .
\end{align*}
From   Kowalczyk et al. \cite{KMM3}  we have
\begin{align}
  - \<  \partial _x^2    \eta ,S_A\eta  \> \ge   2 \|   (\zeta _A \eta)' \| _{L^2}^2- \frac{C }{A}  \|    \eta \| _{\widetilde{\Sigma}} ^2.\label{eq:poscomm}
\end{align}
Obviously
\begin{align}
   |\< V \eta ,S_A\eta  \> |=  2^{-1} |\< [S_A, V] \eta ,\eta  \>|   |\lesssim   \|    \eta \| _{\widetilde{\Sigma}} ^2.\label{eq:potentialcomm}
\end{align}
We write \begin{equation*}
      \<   f(  Q [z ] + \eta  )  -  f(  Q [z ]    )   - f\left(  \eta  \right) ,S_A {\eta} \> + \<  f\left(  \eta  \right) ,S_A {\eta} \>=:B_{ 1}+B_{ 2}.
 \end{equation*}
Then
\begin{align*}
 | B_{ 2}| =  |\<  f\left(  \eta  \right)   \overline{{\eta}} -G (|\eta |^2) , \zeta _A ^2 \> | \lesssim  \int _\R |\eta | ^{p+1} \zeta _A ^2 dx
\end{align*}
and we use   the crucial estimate by    Kowalczyk et al. \cite{KMM3}
\begin{align}\label{eqpurepower}
   \int _\R |\eta | ^{p+1} \zeta _A ^2 dx  \lesssim   A^2 \| \eta \| ^{p-1}_{L^\infty (\R )}   \| (\zeta _A\eta ) ' \| ^{2}_{L^2 (\R )} \lesssim A ^{-1} \| (\zeta _A\eta ) ' \| ^{2}_{L^2 (\R )}.
\end{align}
 We claim
  \begin{align}\label{eq:claim:estB1}
  |B_{ 1}| \lesssim    \|    \eta \| _{\widetilde{\Sigma}} ^2 + \delta _2  \|   (\zeta _A \eta)' \| _{L^2}^2 .
\end{align}  We consider
\begin{align}& \label{eq:claim:estB11}
   \Omega_{1,t  }=\{x\in \R\ |\  |   \eta(t,x)|\leq 2| Q [    z(t)]|\}  \text{ and }\\&  \Omega_{2,t  }=\R\setminus \Omega_{1,t  }=\{x\in \R\ |\ |  \eta(t,x)|> 2| Q [    z(t)]|\} , \nonumber
\end{align}
and split accordingly
\begin{align*}
  B_1 =\sum _{j=1,2} \<   f(  Q [z ]+ \eta  )  -  f(  Q [z ]   )   - f\left(  \eta  \right) ,S_A {\eta} \> _{L^2( \Omega_{j,t  })} =:B _{11}+B _{12}.
\end{align*}
We have
\begin{align*}
  |B _{11}|&\le \int _0 ^1 \int _{\Omega_{1,t  }} \left | \<   \( D   f(   Q [z ] + s\eta  )  -  Df(s \eta   )\)  \eta  ,S_A {\eta}           \> _\C  \right | dx ds\\& \lesssim
  \int _{\Omega_{1,t  }}   \zeta _A ^2| Q [z ]   | ^{p -1} |\eta | ^2 dx +  \int _{\Omega_{1,t  }} \zeta _A ^{-1}  |\varphi  _A |  | Q [z ]  | ^{p -1}  |\eta |    (\zeta _A \eta) ' - \zeta _A '\eta   |      dx =: B _{111}+B _{112}.
\end{align*}
Then
\begin{align*}
  B _{111} \lesssim  \|    \eta \| _{\widetilde{\Sigma}} ^2
\end{align*}
and, for a fixed small and preassigned $\delta _2>0$,
\begin{align*}
  B _{112} \lesssim  \int _{\Omega_{1,t  }}  \zeta _A ^{-1}  |\varphi _A |  |Q [z ]  | ^{p -1}  |\eta |   \(   |(\zeta _A \eta) '|  + \frac{1}{A} \zeta _A |\eta |  \)         dx  \lesssim  \|    \eta \| _{\widetilde{\Sigma}} ^2 + \delta _2  \|   (\zeta _A \eta)' \| _{L^2}^2 .
\end{align*}
We have
\begin{align*}
  |B _{12}|&\le \int _0 ^1 \int _{\Omega_{2,t  }} \left | \<   \( D   f( s Q [z ] + \eta  )  -  Df(s Q [z ]     )\)  Q [z ]   ,S_A {\eta}           \> _\C  \right | dx ds\\& \lesssim
  \int _{\Omega_{2,t  }}   \zeta _A ^2|\eta | ^{p } |Q [z ] | dx +  \int _{\Omega_{2,t  }} \zeta _A ^{-1}  |\phi _A |  |\eta | ^{p-1 } |Q [z ]  |  | (\zeta _A \eta) ' - \zeta _A '\eta   |      dx =: B _{121}+B _{122}.
\end{align*}
We have
\begin{align*}
  B _{121} \lesssim  \int _{\Omega_{2,t  }}   \zeta _A ^2|\eta | ^{p } |Q [z ]  | ^{2-p+p-1}    dx \le  \int _{\R}   |\eta | ^{2 } |Q [z ]  | ^{p-1}    dx \lesssim \|    \eta \| _{\widetilde{\Sigma}} ^2.
\end{align*}
We have similarly, for a fixed small and preassigned $\delta _2>0$ and completing the proof of \eqref{eq:claim:estB1},
\begin{align*}
  B _{122} \lesssim  \int _{\Omega_{2,t  }}  \zeta _A ^{-1}  |\varphi _A | |\eta |   |Q [z ] | ^{p-1}  \(   |(\zeta _A \eta) '|  + \frac{1}{A} \zeta _A |\eta |  \)         dx  \lesssim  \|    \eta \| _{\widetilde{\Sigma}} ^2 + \delta _2  \|   (\zeta _A \eta)' \| _{L^2}^2 .
\end{align*}
 \qed

 \textit{Proof of Proposition \ref{prop:1virial}.}
Integrating  inequality \eqref{lem:1stV1} we obtain \eqref{eq:sec:1virial1}.

\qed

\section{A review of Kato smoothing}\label{sec:smooth10}

This section is mainly inspired by Mizumachi \cite{mizu08} and is based also on material in  \cite{CM2109.08108}   and \cite{CT}.
Since 0 is not a resonance of $H$, by Lemma 1 p. 130 and Theorem 1 \cite{DT}  and \cite[formula (2.45)]{weder} we have the following result about \textit{Jost functions}.
  \begin{proposition}[Jost functions]
   \label{prop:DT} For any $k\in \overline{ \C }_+ $ there exist functions   $f_{\pm } (x,k )=e^{\pm \im k x}m_{\pm } (x,k )$  which solve $ H u=k^2 u$ with\begin{align*} \lim _{x\to +\infty }   {m_{ + } (x,k )}  =1 =\lim _{x\to -\infty }  {m_{- } (x,k)}  .\end{align*}  These functions, for $x^{\pm }:= \max \{ 0,\pm  x \}$  satisfy\begin{align}  \label{eq:kernel2} &  |m_\pm(x, k )-1|\le  C _1 \langle  x^{\mp}\rangle\langle k \rangle ^{-1}\left | \int _x^{\pm \infty}\langle y \rangle |V  (y)| dy \right |    .   \end{align}
   For the Wronskian we have $     [f_{ +} (x,k ), f_{ -} (x,k ) ] =   \dfrac{2\im k}{T(k)}  $ where  $T(k) =\alpha k (1+o(1))$  near $k=0$ for some $\alpha \in \R$     and $T(k) = 1+O(1/k) $ for $k\to \infty$ and $T\in C^0(\R )$.
\end{proposition}
    \qed

  For $k\in  { \C }_+ $  we   have   the following formula for the integral the kernel of the resolvent $R_{H}  (k^2)$ of $H$,                  for $x<y$, with an analogous formula for $x>y$,
   \begin{align}&   R _{H }    (k^2 ) (x,y) = \left\{
                                                  \begin{array}{ll}
                                                     \frac{T(k)}{2\im k}     f_- (x,k)   f_+ (y, k)   =   \frac{T(k)}{2\im k}     e^{\im k (x-y)}       m_- (x, k)   m_+ (y, k)  & \hbox{ for $x<y$ } \\
                                                     \frac{T(k)}{2\im k}     f_+ (x,k)   f_+ (-, k)   =   \frac{T(k)}{2\im k}      e^{-\im k (x-y)}       m_+ (x, k)   m_- (y, k)  & \hbox{ for $x>y$ .}
                                                  \end{array}
                                                \right.
   \label{eq:resolvKern} \end{align}
Let now  $\lambda \ge 0$ and set $k= \sqrt{\lambda}\ge 0$.
Then the   formula     \eqref{eq:resolvKern}   makes sense   yielding the kernel $R _{H }  ^+  (\lambda  ) (x,y)$  of an operator which is denoted by $R_{H} ^+  (\lambda )$.
We can define similarly Jost functions for $\Im k<0$  which are of the form $\overline{f_{\pm } (x,\overline{k} )}$. In particular $ R _{H } ^{-}   (k^2 ) (x,y)= \overline{R} _{H } ^{+}   (k^2 ) (x,y)$  for $k\ge 0$.

 We have the following.

\begin{lemma}
  \label{lem:5.7} For any $s>3/2$   and $\tau >1/2$ there is a constant $C _{s\tau }>0$ such that for  any $\lambda \ge 0$  the following limit exists for
\begin{align}
   \label{eq:5.71}   & \lim _{a \to 0 ^+} R_{H }(\lambda \pm  \im a )
=R^{\pm }_{H }(\lambda  ) \text{ in }   \mathcal{L}  \(  L ^{2,\tau  } \( \R \)  , L ^{2,-s } \( \R \)   \)   \text{  and we have the bound }
\\&  \sup _{\lambda , a \in \R _+} \|    R_{H }(\lambda \pm  \im a ) \| _{\mathcal{L}  \(  L ^{2,\tau } \( \R \)  , L ^{2,-s } \( \R \)   \)} \le C _{s\tau } . \label{eq:5.72}
\end{align}
\end{lemma}
 \proof    The result is standard and we sketch it for completeness. It is enough to consider the $+$ sign.
For any $\lambda \ge 0$,  $a>0$ small,  $x<y$       and   $k_{a} = \sqrt{\lambda + \im a }$ we have  \small
\begin{align} \label{eq:5.73}
   \< x \> ^{-s}\left |  R _{H }    (\lambda + \im a ) (x,y)  \right |  \< y \> ^{-\tau} \lesssim   \< x \> ^{-s}     \(  1+ x^+ + y^-  \)  \< y \> ^{-s}\lesssim  \left\{
                                                                                                                                                                  \begin{array}{ll}
                                                                                                                                                                      \< x \> ^{-s+1}  \< y \> ^{-\tau}  & \hbox{ if $y>x>0$;} \\
                                                                                                                                                                       \< x \> ^{-s }  \< y \> ^{-\tau+1} & \hbox{ if $x<y<0$.}\\
                                                                                                                                                                      \< x \> ^{-s }  \< y \> ^{-\tau} & \hbox{ otherwise.}
                                                                                                                                                                  \end{array}
                                                                                                                                                                \right.
\end{align}\normalsize
We have an analogous estimate for $x>y$.
It is then  elementary  to show that the latter kernel   $\< x \> ^{-s}   R _{H }    (\lambda + \im a ) (x,y)     \< y \> ^{-\tau } \in L^{2}\(   \R \times \R \)$ with norm independent from $a>0$ and $\lambda \ge 0$.
Indeed
for $z=\lambda \pm \im a$  it is enough to  consider
 \begin{align} \nonumber &   \int _{\R} dx \< x \> ^{-2s} \int_{\R}  |R  _{H }(x,y,z )| ^2  \< y \> ^{-2\tau} dy   =    \int _{\R} dx \< x \> ^{-2s} \int_{-\infty}^{x}  |R  _{H }(x,y,z )| ^2  \< y \> ^{-2\tau} dy\\& + \int _{\R} dx \< x \> ^{-2S} \int_{x}^{+\infty}  |R  _{H }(x,y,z )| ^2  \< y \> ^{-2\tau} dy.\label{eq:estresolv11}
\end{align}
 The second term in the right hand side  is bounded   using  \eqref{eq:5.73}   by
  \begin{align*}  &       \int _{x<y} \< x \> ^{-2s}  \< y \> ^{-2\tau}  \( 1+ x^++  y^- \) ^2   dx  dy
   \le  \int _{0<x<y } \< x \> ^{-2s+2}  \< y \> ^{-2\tau  }    dx  dy   \\&+ \int _{x<y<0} \< x \> ^{-2s}  \< y \> ^{-2\tau +2}    dx  dy + \int _{ x<0<y } \< x \> ^{-2s }  \< y \> ^{-2\tau  }    dx  dy =:\sum _{j=1}^{3}I_j.
\end{align*}
Then
\begin{align*}
  I_1\le  \int _{\R } \< x \> ^{-2S+2}  dx \int _\R   \< y \> ^{-2\tau  }      dy =:I_4 <\infty \text{  for $s>3/2$ and $\tau >1/2$.}
\end{align*}
Similarly $I_j< I_4$ for $j=2,3$. Similar estimates  independent from $\lambda >0$ and $a>0$  hold for the term in the first line in the right hand side of \eqref{eq:estresolv11}. This gives us  also \eqref{eq:5.72}. The exact same bound   holds  for the $L^2(\R \times \R )$ norm of
 $\< x \> ^{-s}  R _{H }  ^+   (\lambda   ) (x,y)    \< y \> ^{-\tau}$.   Furthermore
\begin{align*}
   \lim _{a \to 0 ^+}  \< x \> ^{-s}   R _{H }    (\lambda + \im a ) (x,y)     \< y \> ^{-\tau} = \< x \> ^{-s}  R _{H }  ^+   (\lambda   ) (x,y)    \< y \> ^{-\tau} \text{ for all $x,y$}
\end{align*}
and by Lebesgue's  dominated convergence the above convergence holds in $L^{2}\(   \R \times \R \)$. This yields the limit \eqref{eq:5.71}  for any $\lambda \ge 0$.
\qed

Taking $s=\tau$ we  have the following.
\begin{lemma}[Kato smoothing]
  \label{lem:3.3}  For  any $s>3/2 $  there exists a constant  $c_s>0$ such that
\begin{align*}
  \| e^{-\im tH }P_c f\|
_{  L^2 \(\R  ,  L ^{2 , -s } (\R )\) } \le
c_s    \|f\|_{ L ^{2   } (\R )};
\end{align*}
\end{lemma}
\proof  We sketch this well known result for completeness.
Get $g (t,x) \in \mathcal{S}(  \R \times \R ) $ with $g(t)=P_c g(t)$. By the limiting
absorption principle, taking $f$ such that $f=P_cf$
\begin{align*}
  \(  e^{-\im tH  } f,
g\)  _{L^2(\R _t\times \R _x )}  =    \frac{1}{ 2\pi \im }\int_0 ^{+\infty}
   \(
 (R_{H  }^{+}(\lambda  )-R_{H
}^{-}(\lambda  ))   f,     \widehat{g} (\lambda
)\) _{L^2 (\R _x)} d\lambda .
\end{align*}
Then from Fubini and  Plancherel   we have
\begin{align} \nonumber
   \left |\( e^{-\im tH  } f, g\) _{L^2(\R _t\times \R _x ) }\right  |
 &\le (2\pi)^{-\frac{1}{2}}
\| (R_{H  }^{+}(\lambda  )-R_{H  }^{-}(\lambda )) f\|_{     L ^{2,
-s} _x L^2_\lambda   }   \|  {g} \|_{L ^{2,
 s} _x  L^2_t } \\&= (2\pi)^{-\frac{1}{2}}\| (R_{H  }^{+}(\lambda  )-R_{H  }^{-}(\lambda )) f\|_{ L^2\( \R ,     L ^{2,
-s} \( \R \)\)    }   \|  {g} \|_{L^2\( \R ,     L ^{2,
 s} \( \R \)\)   }
 .\label{eq:3.2}
\end{align}
By Fathou lemma
\begin{align*}
  \| (R_{H  }^{+}(\lambda  )-R_{H  }^{-}(\lambda )) f\|_{L^2\( \R ,     L ^{2,
-s} \( \R \)\)    }\le \liminf _{a \to 0 ^+} \| (R_{H }
(\lambda +\im   a )-R_{H  } (\lambda -\im   a  )) f\|_{L^2\( \R ,     L ^{2,
-s} \( \R \)\)  }.
\end{align*}
 Then we repeat an argument in  \cite[Lemma 5.5]{kato}. For $\Im \mu >0$
set $K(\mu )$ the positive square root of $(2\pi \im )^{-1} [R_{H  }
(\mu )- R_{H } ( \overline{\mu} )]= \pi ^{-1}(\Im \mu )R_{H  }
(\overline{\mu} ) R_{H  } ( {\mu} )$. Then \small
\begin{align*}&  \int  _{\R _+} \|(R_{H  }^{+}(\lambda  )-R_{H  }^{-}(\lambda )) f\|_{ L ^{2,
-s} \( \R \)  } ^2d\lambda  \le
   \liminf _{a \to 0 ^+} \int  _{\R _+}    \| (R_{H }
(\lambda +\im   a )-R_{H  } (\lambda -\im   a)) f\|_{ L ^{2,
-s} \( \R \)  } ^2d\lambda  \\& =4\pi ^2 \liminf _{a \to 0 ^+} \int _{\R _+}    \| \langle x \rangle
^{-s }K(\lambda +\im  a )K (\lambda +\im  a  )
f\|_{L ^{2 } \( \R \)  } ^2d\lambda
  \\& =  4\pi ^2 \liminf _{a \to 0 ^+} \int _{\R _+}    \|  K(\lambda +\im  a )    \langle x \rangle
^{-s }  K (\lambda +\im  a  )
f\|_{L ^{2 } \( \R \)  } ^2d\lambda \\& = -2\im \pi  \liminf _{a \to 0 ^+}   \int _{\R _+}   \(   \langle x \rangle
^{-s } \(  R_{H }
(\lambda +\im   a )-R_{H  } (\lambda -\im   a)  \)   \langle x \rangle
^{-s }  K (\lambda +\im  a  )
f ,        K (\lambda +\im  a  )
f \)  d\lambda \\& \le 4  \pi C_{ss}  \liminf _{a \to 0 ^+}  \int _{\R _+}      \|  K
(\lambda +\im a ) f\|_{L ^{2 } \( \R \) } ^2d\lambda  \\& =  4  \pi C_{ss}    \(  \lim _{a \to 0 ^+} \frac{1}{2\pi \im }  \int _{\R _+}  \(  R_{H }
(\lambda +\im   a )-R_{H  } (\lambda -\im   a)   \)  fd\lambda , f\)  =  4  \pi C_{ss}    \| f\|_{L ^{2 } \( \R \)  } ^2 .
\end{align*}
\normalsize
Then by    \eqref{eq:3.2} we get the following, which yields the lemma,
\begin{align*}
   &  \left |\( e^{-\im tH  } f, g\) _{L^2(\R _t\times \R _x ) }\right  | \le   c_s    \| f\|_{L ^{2 } \( \R \)  }  \|  {g} \|_{L^2\( \R ,     L ^{2,
 s} \( \R \)\)   }  .
\end{align*}
 \qed

The following formulas are taken from Mizumachi \cite[Lemma 4.5]{mizu08}   to which we refer for the proof.
 \begin{lemma} \label{lem:lemma11} Let  for $g\in \mathcal{S}(\R \times \R , \C )$   with $g(t) \equiv P_cg(t)$,  \begin{align*}&  U(t,x) = \frac{1}{\sqrt{2\pi} \im }\int _\R e^{-\im \lambda t}\( R ^{-}_{H }(\lambda )+R ^{+}_{H_{N+1}}(\lambda )   \)  \breve{g}( \lambda   ) d\lambda  ,
\end{align*}
where $\breve{g}( \lambda   )$  is the inverse Fourier transform in $g(t)$ (in $t$). Then
 \begin{align} \label{eq:lemma11}2\int _0^t e^{-\im (t-t') H  }g(t') dt'   &=  U(t,x)  - \int  _{\R _-} e^{-\im (t-t') H  }g(t') dt'
 \\& + \int  _{\R _+} e^{-\im (t-t') H  }g(t') dt' .\nonumber
\end{align}
\end{lemma}
\qed

The following smoothing is due to Mizumachi \cite{mizu08}.
\begin{lemma}\label{lem:smooth}  For  any $s>3/2 $  and $\tau >1/2$ there exists a constant  $ c_{s\tau }>0$ such that
   \begin{align}\label{eq:smoothest1}
      \left\|  \int_0^t e^{-\im (t-s)H }
P_c(H)g(s)ds\right\|_{L^2\( \R ,     L ^{2,
 -s} \( \R \)\)} \le c_{s\tau }   \|
g\|_{L^2\( \R ,     L ^{2,
 \tau } \( \R \)\)}.
   \end{align}
\end{lemma}
\proof For completeness we sketch the proof, which is similar to  \cite[Lemma 8.9]{CM24D1}.       We can use formula \eqref{eq:lemma11} and bound $U$, with the bound on the last two terms in the right hand side of \eqref{eq:lemma11} similar. Taking Fourier transform in $t$,
$$ \aligned &
\| U\|_{L^2\( \R ,     L ^{2,
 -s} \( \R \)\)} \le   2 \sup _{\pm}\| R_{ H  }^\pm (\lambda )
   \widehat{ g}(\lambda  )\|_{L^2\( \R ,     L ^{2,
 -s} \( \R \)\) }   \le \\& \le
2\sup _{\pm}   \sup _{\lambda \in \R }
\|  R_{ H  }^\pm (\lambda )  \| _{ L^{2,\tau}  \to L^{2,-s}   } \|
     \widehat{g} (\lambda ) \|_{ L^2\( \R ,     L ^{2,
 \tau} \( \R \)\) }\,
 \lesssim  &        \| g\|_{L^2\( \R ,     L ^{2,
 \tau} \( \R \)\) }.
\endaligned $$
Notice that, while Lemma \ref{lem:lemma11} is stated for $g\in \mathcal{S}(\R \times \R , \C )$, the estimate \eqref{eq:smoothest1} extends to all $g\in L^2( \R , L^{2,\tau}(\R ) )$ by density.
\qed

\section{Low energies: proof of Proposition \ref{prop:smooth11}}\label{sec:smooth1}

Setting   $v:= \chi _B   \eta $  for $B$ like in \eqref{eq:relABg}  and then $w=P_c v$ from \eqref{eq:nls3}    we have

\begin{align}       &
    \im \dot w   =      H   w +P_c \( \chi _B'' +2\chi _B'\partial _x  \) \eta    -  \im  P_c     \chi _B   D_{z} Q [z ]  \(\dot  z   +\im E z \)  +   P_c     \chi _B   \( f(  Q [z ]  + \eta  )  -  f(  Q [z ]  )  \)     .   \label{eq:nls5}
\end{align}
Then we have the following, using the spaces in \eqref{eq:l2w}.

\begin{lemma}\label{lem:estw}
 For $s > 3/2$
   we have
 \begin{equation}\label{eq:estw1}
   \| w \|  _{L^2(I , L ^{2,-s}(\R ) )}\le   o_{B^{-1}}(1)   \epsilon.
 \end{equation}

\end{lemma}  \proof
By  Lemma \ref{lem:3.3} for the first inequality  and using $w=P_cw$,        we have
 \begin{align}
   \|    e ^{\im t H }   w(0) \| _{L^2(\R, L ^{2,-s}(\R ))} \lesssim   \|
   w(0)  \| _{L^2 (\R )} \lesssim \|
   \eta(0)  \| _{L^2 (\R )} \lesssim \delta \( = o_{B^{-1}}(1)   \epsilon  \).\nonumber
 \end{align}
 Using Proposition \ref{prop:bddst} and   Lemmas \ref{lem:lemdscrt} and \ref{lem:smooth}, taking $\tau\in (1/2,1)$ we get
\begin{align*}&
   \|  \int _{0}^t e^{-\im (t-t') H  }P_c     \chi _B   D_{z} Q [z ]  \(\dot  z   +\im E z \)    dt' \| _{L^2(I , L ^{2,-s}(\R ) )}  \lesssim    \|     \dot  z   +\im E z   \| _{L^2(I   )}     \|     D_{z} Q [z ]    \| _{L^\infty (I , L ^{2, \tau }(\R ) )} \\& \lesssim    \delta ^{p-1}   \|       \eta   \| _{L^2\(I  , \widetilde{\Sigma} \)}   = o_{B^{-1}}(1)   \epsilon .
\end{align*}
 By   Lemma  \ref{lem:smooth}
\begin{align*}&
   \|  \int _{0}^t e^{-\im (t-t') H} P_c \( \chi _B'' +2\chi _B'\partial _x  \) \eta  dt'    \| _{L^2(I , L ^{2,-s}(\R ) )} \\& \lesssim
  \|     \( 2 \chi ' _B \partial _x + \chi _B '' \)  \eta     \| _{L^2(I , L ^{2,\tau }(\R ) )}    \lesssim B ^{\tau-1} \| \sech \( \frac{2}{A}x\) \eta ' \| _{L^2(I ,L^2(\R ))} \\& + B ^{\tau -2} \| 1_{B\le |x|\le 2 B} \sech \( \frac{2}{A} x\)  \eta   \| _{L^2(I ,L^2(\R ))} \lesssim B ^{\tau-1} \|   \eta   \| _{L^2(I ,\Sigma _A )} \\& +  B ^{\tau -1} \( \left \|   \( \sech \( \frac{2}{A}x\) \eta  \) '    \right \| _{L^2(I ,L^2(\R ))}   +  \|    \eta \| _{L^2(I , \widetilde{\Sigma})}
  \) =o_{B^{-1}} (1)   \epsilon,
\end{align*}
where we used   $\tau \in (1/2,1) $   and, see   Merle and Raphael \cite[Appendix C]{MR4},
\begin{align*}
  \|   u \| _{L^2(|x|\le 2 B)}\lesssim  B \( \left \|  u '    \right \| _{ L^2(\R  )}   +  \|    u \| _{  \widetilde{\Sigma} }\)  .
\end{align*}
 Next we consider
\begin{align*} &
    \|  \int _{0}^t e^{-\im (t-t') H} P_c     \chi _B   \( f(  Q [z ]  + \eta  )  -  f(  Q [z ]  \)   dt'    \| _{L^2(I , L ^{2,-s}(\R ) )} \le I_1+I_2 \text{ where} \\& I_1=  \|\int _{0}^t e^{-\im (t-t') H} P_c     \chi _B     f(   \eta  )   dt'    \| _{L^2(I , L ^{2,-s}(\R ) )}   \text{ and} \\& I_2=  \|\int _{0}^t e^{-\im (t-t') H} P_c     \chi _B   \( f(  Q [z ]  + \eta  )  -  f(  Q [z ]  ) -f(   \eta  ) \)      dt'    \| _{L^2(I , L ^{2,-s}(\R ) )} .
\end{align*}
We have by \eqref{eq:relABg}
\begin{align*}
  I_1&\lesssim   \|      \chi _B     f(   \eta  )       \| _{L^2(I , L ^{2,\tau }(\R ) )}  \lesssim   B ^{\tau}\| \eta \|  ^{p-1} _{L^\infty(I, L^2 (\R ) )} A    A ^{-1} \|     \sech \( \frac{2}{A} x\)  \eta    \| _{L^2(I , L ^{2  }(\R ) )} \\& \lesssim \delta ^{p-1}  B ^{\tau}  A    \|        \eta    \| _{L^2(I , \Sigma _A )} = o_{B^{-1}} (1)   \epsilon .
\end{align*}
We have
\begin{align*}
  I_2 ^2 & \lesssim  B ^{2\tau}   \int _{I\times \R }    \chi _B^2    \left |  \int _0 ^1 \left [ Df(  Q [z ]  + \sigma  \eta  )    -Df( \sigma  \eta  )  \right ] \eta \right  |  ^2 dx dt   .
\end{align*}
Then by  $|Df(  Q [z ]  + \sigma  \eta  )    -Df( \sigma  \eta  )|  \le  |Df(  Q [z ]  + \sigma  \eta  ) |   + |Df( \sigma  \eta  )| \lesssim   |Q[z]| ^{p-1}   +|  \eta |^{p-1}$,
\begin{align*}
   I_2 ^2 & \lesssim  B ^{2\tau}  \int _{I\times \R }    \chi _B^2    \( |z| ^{2p-2}   + |  \eta    | ^{2p-2} \) |\eta |  ^2  dx dt   \\& \lesssim   B ^{2\tau}
\( |z| ^{2p-2}   +\|  \eta  \| ^{2p-2}_{L^\infty ( I,   L^\infty ( \R ) )} \)  A A ^{-1}      \| \sech \( \frac{2x}{A}  \)\eta \| _{L^2 ( I,   L^2 ( \R ) )}\lesssim  B ^{2\tau} A \delta ^{p-1}  \| \eta \| _{L^2 ( I,   \Sigma_{A} )}.
\end{align*}
So we conclude, the following, which completes the proof of   \eqref{eq:estw1},
\begin{align*}
  I_1+I_2  = o_{B^{-1}} (1)   \epsilon  .
\end{align*}

\qed

\textit{Proof of Proposition \ref{prop:smooth11}.}
We have  by $P \eta =0$
\begin{align*}
  v= \chi _B \eta  = P_c v + P  \chi _B \eta = w  - \varphi \<  (1- \chi _B) \varphi , \eta \>   - \im \varphi \< \im  (1- \chi _B) \varphi , \eta \>
\end{align*}
  where   $ \varphi (x) \lesssim e^{-\sqrt{\lambda}|x|} $ yields the pointwise  estimate
\begin{align*}
  |\varphi \<  (1- \chi _B) \varphi , \eta \>   + \im \varphi \< \im  (1- \chi _B) \varphi , \eta \>| \lesssim e^{-\sqrt{\lambda}B}     \| \eta \| _{   \widetilde{\Sigma} } .
\end{align*}
This and \eqref{eq:estw1} yield   $\| v\| _{L^2(I , \widetilde{\Sigma})}\lesssim o_{B^{-1}}(1) \epsilon$.
Next, from  $v = \chi _B    \eta $,  and thanks to the relation  $A\sim B^3$  set in \eqref{eq:relABg}    we have
\begin{align} \nonumber
   \| \eta \| _{\widetilde{\Sigma}}&\lesssim  \| v \| _{\widetilde{\Sigma}} + \| (1- \chi _B   )\eta \| _{\widetilde{\Sigma}}\lesssim \| v \| _{\widetilde{\Sigma}} + A ^{-2} \|  \sech \( \frac{2}{A}x  \) \eta \| _{L^2}\\& \lesssim \| v \| _{\widetilde{\Sigma}} + A ^{-1} \|   \eta \| _{\Sigma _A}\label{eq:relAB0}
\end{align}
So by \eqref{eq:sec:1virial1}    we get  the following, which implies \eqref{eq:sec:smooth11},
\begin{align} \nonumber
   \| \eta \| _{L^2(I,\widetilde{\Sigma})}&\lesssim    \| v \| _{L^2(I,\widetilde{\Sigma})} + A ^{-1} \|   \eta \| _{L^2(I, \Sigma _A)}\\& \lesssim  o_{B^{-1}} (1) \epsilon +  A ^{-1} \( \| \dot  z   +\im E z\|_{L^2(I)} + \| \eta \| _{L^2(I, \widetilde{\Sigma}   )}\) \lesssim   o_{B^{-1}}(1) \epsilon +  A ^{-1}\| \eta \| _{L^2(I, \widetilde{\Sigma}   )}    .  \nonumber
\end{align}

 \qed

\section*{Acknowledgments}
C. was supported   by the Prin 2020 project \textit{Hamiltonian and Dispersive PDEs} N. 2020XB3EFL.
M.  was supported by the JSPS KAKENHI Grant Number 19K03579, 23H01079 and 24K06792.

Department of Mathematics, Informatics and Geosciences,  University
of Trieste, via Valerio  12/1  Trieste, 34127  Italy.
{\it E-mail Address}: {\tt scuccagna@units.it}

Department of Mathematics and Informatics,
Graduate School of Science,
Chiba University,
Chiba 263-8522, Japan.
{\it E-mail Address}: {\tt maeda@math.s.chiba-u.ac.jp}

\end{document}